\documentclass[11pt]{amsart}
\usepackage{amsmath,amsfonts,amssymb}
\usepackage{url}

\newtheorem{theorem}{\sc Theorem}

\newcommand{\group}{\mathcal{G}}
\newcommand{\graph}{\mathbf{G}}
\newcommand{\cantor}{2^{\mathbb{N}}}
\newcommand{\Aut}{\text{Aut}}
\newcommand{\Homeo}{\text{Homeo}}
\newcommand{\hilbert}{[0,1]^{\mathbb{N}}}

\usepackage{tikz}
\usepackage{tikz-cd}

\begin{document}
\title{Frucht's theorem in Borel setting}
\author{Onur B\.{i}lge}
 \address{Department of Mathematics \\ Gazi University, 06500, Ankara, Turkey}
 \email{onurbilge@gazi.edu.tr}
\address{Department of Mathematics \\ Middle East Technical University,
 06800, Ankara, Turkey}
\email{bilge.onur@metu.edu.tr}
\author{Burak KAYA}
\address{Department of Mathematics \\ Middle East Technical University,
 06800, Ankara, Turkey }
\email{burakk@metu.edu.tr}

\begin{abstract}
In this paper, we show that Frucht's theorem holds in Borel setting. More specifically, we prove that any standard Borel group can be realized as the Borel automorphism group of a Borel graph. A slight modification of our construction also yields the following result in topological setting: Any Polish group can be realized as the homeomorphic automorphism group of a $\mathbf{\Delta^0_2}$-graph on a Polish space.
\end{abstract}
\maketitle

\section{Introduction}
Descriptive graph combinatorics is the study of ``definable" graphs on Polish spaces while incorporating descriptive-set-theoretic concepts into graph-theoretic concepts. Working on \textit{a Polish space} i.e. a completely metrizable separable topological space, one can require various graph-theoretic objects such as edge relations, colorings, perfect matchings, automorphisms etc. to have topological and measure-theoretic properties such as being Borel, projective, continuous, closed etc. One can then ask to what extent classical results from graph theory in abstract setting generalize to measurable setting. Some classical theorems, e.g. the $(k+1)$-colorability of a locally finite graph of bounded degree $k$, generalizes to definable setting \cite[Proposition 4.6]{KST99} while some others, e.g. the $2$-colorability of an acyclic graph, do not \cite[3.1. Example]{KST99}. It turns out that this jump from the abstract setting to the measurable setting is more than a mere specialization and leads to interesting and fruitful results. We refer the reader to \cite{KechrisMarks20} for a comprehensive treatment of the subject.

This study has so far focused on chromatic numbers and perfect matchings. Analyzing measurable automorphisms groups of graphs seems to be a natural direction to extend this study. In this paper, we shall tackle the problem of realizing a measurable group as the measurable automorphism group of a definable graph.

A classical theorem of Frucht in \cite{Frucht39} states that every finite group is isomorphic to the automorphism group of some finite graph. In \cite{Groot59} and \cite{Sabidussi60}, de Groot and Sabidussi independently generalized Frucht's theorem to arbitrary groups by removing the finiteness condition. All these result seem to use Frucht's original idea that can be summarized as follows:

\begin{quote}
        Given a group $\group$ with a generating set $S$, consider the Cayley graph $G$ with respect to $S$ as a directed labeled graph. Then the group of automorphisms of $G$ as a directed labeled graph is isomorphic to $\group$. Systematically replace each directed labeled edge by a connected undirected asymmetric graph to obtain an undirected graph. Then the automorphism group of the resulting undirected graph is isomorphic to $\group$.
\end{quote}

While this idea does not seem to invoke any non-explicit methods at first glance, such as the use of the axiom of choice that often results in non-measurable objects, it remains a non-trivial question to answer whether or not the ``systematically replace" part of this idea can actually be done in a uniform way in Borel setting. Indeed, the arguments in \cite{Groot59} and \cite{Sabidussi60} do not seem to produce Borel graphs. Nevertheless, the answer turns out to be affirmative as we shall see later. Before we give a precise statement of our main result, let us recall some basic notions from descriptive set theory. We refer the reader to \cite{Kechris95}  for a general reference.

A measurable space $(X,\mathcal{B})$ is called a \textit{standard Borel space} if $\mathcal{B}$ is the Borel $\sigma$-algebra of a Polish topology on $X$. A graph $\graph=(X,G)$ on a standard Borel space $(X,\mathcal{B})$ is said to be \textit{Borel} if its edge relation $G \subseteq X \times X$ is a Borel subset of the product space. The group of automorphisms of $\graph$ that are Borel maps will be denoted by $\Aut_B(\graph)$. In the case that $\graph$ is a Borel graph on a Polish space $(X,\tau)$, the group of automorphisms of $\graph$ that are homeomorphisms will be denoted by $\Aut_h(\graph)$.

A triple $(\group,\cdot,\mathcal{B})$ is said to be a \textit{standard Borel group} if $(\group,\mathcal{B})$ is a standard Borel space and $(\group,\cdot)$ is a group for which the multiplication $\cdot: \group \times \group \rightarrow \group$ and the inversion $\ ^{-1}: \group \rightarrow \group$ operations are Borel maps. A triple $(\group,\cdot,\tau)$ is said to be a \textit{Polish group} if $(\group,\tau)$ is a Polish space and $(\group,\cdot)$ is a group for which $\cdot: \group \times \group \rightarrow \group$ and $\ ^{-1}: \group \rightarrow \group$ are continuous maps. The main result of this paper is the following variation of Frucht's theorem in Borel setting.

\begin{theorem}\label{maintheorem} For every standard Borel group $(\group,\cdot,\mathcal{B})$, there exists a Borel graph $\graph=(X,G)$ on a standard Borel space $(X,\widehat{\mathcal{B}})$ such that $\group$ and $\Aut_B(\graph)$ are isomorphic.
\end{theorem}

A slight modification of our argument in the proof of Theorem \ref{maintheorem} also gives the following variation in topological setting.

\begin{theorem}\label{maintheorem2} For every Polish group $(\group,\tau)$, there exists a $\mathbf{\Delta^0_2}$-graph $\graph=(X,G)$ on a Polish space $(X,\widehat{\tau})$ such that $\group$ and $\Aut_h(\graph)$ are isomorphic. Moreover, this isomorphism can be taken to be a homeomorphism where $\Aut_h(\graph) \subseteq \Homeo(X)$ is endowed with the subspace topology induced from the compact-open topology of $\Homeo(X)$.
\end{theorem}

As it was hinted before, our construction is a Borel implementation of Frucht's original idea with appropriate coding techniques which makes sure that the resulting edge relation stays Borel. This paper is organized as follows. In Section 2, after supplying continuum-many pairwise non-isomorphic asymmetric connected countable graphs to be used as replacements of labeled directed edges, we shall construct our candidate graph and prove that its edge relation is Borel. In Section 3, we shall prove Theorem \ref{maintheorem}. In Section 4, recasting the proof of Theorem \ref{maintheorem} with appropriate modifications, we will prove Theorem \ref{maintheorem2}. In Section 5, we shall discuss some further directions and open questions on this theme that can be explored.

\textbf{Acknowledgements.} This paper is a part of the first author’s master’s
thesis \cite{Bilge22} written under the supervision of the second author at the Middle East
Technical University.

\section{Constructing the graph}

Throughout the paper, $\cantor$ denotes the Cantor space i.e. the Polish space consisting of binary sequences indexed by natural numbers, $R^*$ denotes the symmetrization of a relation $R$ on a set i.e. $R^* = R \cup R^{-1}$, $\Delta_X$ denotes the identity relation on $X$ and $\mathbb{N}_{\geq k}$ denotes the set of natural numbers greater than or equal to $k$.

For each $\mathbf{a} \in \cantor$, consider the graph $\graph_{\mathbf{a}}=(\mathbb{N}_{\geq 2},R_{\mathbf{a}}^*)$ where the edge relation is the symmetrization of the relation $R_{\mathbf{a}}=A_{\text{initial}} \cup A_{\text{fork}} \cup A_{\text{nofork}}$ with
\begin{align*}
A_{\text{initial}}=&\left\{(2,3),(3,4)\right\}\\
A_{\text{fork}}=&\left\{(n,n+1),(n,n+2):\ n \in 2\mathbb{N}_{\geq 2},\ \mathbf{a}\left(\frac{n-4}{2}\right)=1\right\}\\
A_{\text{nofork}}=&\left\{(n,n+1),(n+1,n+2):\ n \in 2\mathbb{N}_{\geq 2},\ \mathbf{a}\left(\frac{n-4}{2}\right)=0 \right\}
\end{align*}

The placement of edges in $\graph_{\mathbf{a}}$ can be described as an iterative process as follows. Regardless of $\mathbf{a}$, we first put an edge between $2$ and $3$, and, $3$ and $4$. For each even integer $n \geq 4$, depending on whether $\textbf{a}\left(\frac{n-4}{2}\right)$ is zero or one, we either create a fork at $n$ using the next two vertices with odd vertex having degree one, or add an edge between successive vertices for the next two vertices. For example, a diagrammatic representation of $\graph_{\mathbf{a}}$ with   $\mathbf{a}=(1,0,1,1,0,\dots)$ is as follows.\\

\begin{center}
    \includegraphics[width=0.9\textwidth]{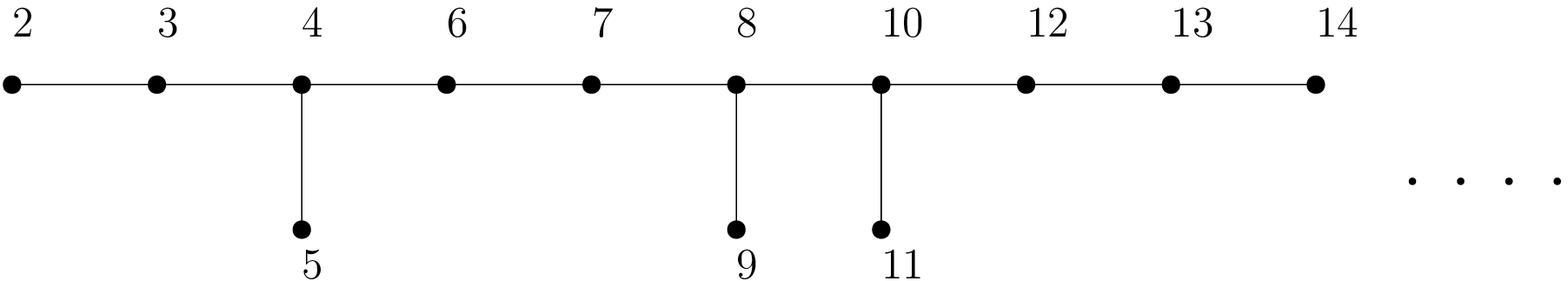}
\end{center}

We shall now argue that any such graph $\graph_{\mathbf{a}}$ is asymmetric i.e. it has no non-trivial automorphisms. Let $\mathbf{a} \in \cantor$ and $\varphi \in Aut\left(\graph_{\mathbf{a}}\right)$. Observe that $2$ is the only vertex of degree one that is adjacent to a vertex of degree two. Hence $\varphi$ fixes $2$ which immediately implies that $3$ and $4$ are fixed under $\varphi$ as well. Let $n \geq 4$ be an even integer. Suppose that $\varphi$ fixes all vertices $2 \leq k \leq n$. Then there are two possibilities:
\begin{itemize}
    \item If $\textbf{a}\left(\frac{n-4}{2}\right)=1$, then $n+1$ is a vertex of degree one and $n+2$ is a vertex of degree two, in which case $\varphi$ fixes both.
    \item If $\textbf{a}\left(\frac{n-4}{2}\right)=0$, then $\varphi$ clearly fixes $n+1$ because $n$ is fixed by $\varphi$ and the other neighbors of $n$ have already been fixed. But subsequently, $\varphi$ must fix $n+2$ as well by a similar argument.
\end{itemize}
Therefore $\varphi$ fixes all the vertices $2 \leq k \leq n+2$. By induction, $\varphi$ fixes all vertices in $\graph_{\mathbf{a}}$. A similar inductive argument shows that $\graph_{\mathbf{a}}$ and $\graph_{\mathbf{b}}$ are not isomorphic whenever $\mathbf{a}$ and $\mathbf{b}$ are distinct elements of $\cantor$.

Next will be constructed the main graph associated to an uncountable standard Borel group. Fix an uncountable standard Borel group $(\group,\cdot,\mathcal{B})$. In order to implement Frucht's idea, we first need to find an appropriate Cayley graph for $(\group,\cdot,\mathcal{B})$. An obvious choice for a generating set is the Borel set $S=\group\setminus\{1_{\group}\}$. Suppose that we constructed the Cayley graph associated to this generating set. In this graph, there is a labeled directed edge from the first component to the second components of each element of $(\group \times \group) \setminus \Delta_{\group}$. We would like to replace each of these directed labeled edges by an appropriate asymmetric connected countable graph that we have already constructed. Consequently, for each element of $(\group \times \group) \setminus \Delta_{\group}$, we need to add countably many ``new" vertices to ``old" vertices. Therefore, it is natural to consider
\[X=\group\times\group\times\mathbb{N}\]
as the vertex set of the main (undirected) graph to be constructed. In this vertex set,
\begin{itemize}
    \item  the vertices of the form $(x,x,0)$ where $x\in \group$ are supposed to represent the ``old" vertices that are the group elements,
    \item the vertices of the form $(x,y,k)$ where $x \neq y \in \group$ and $k \in \mathbb{N}$ are the ``new" vertices that are added after replacing the directed labeled edges, and
    \item the vertices of the form $(x,x,k)$ where $x \in \group$ and $k \neq 0$ are ``irrelevant" elements that will essentially serve no purpose. We could simply have taken these elements out of the vertex set, however, there is no harm in keeping them around. In order for these vertices to not create any additional symmetries, we will stick an infinite line formed by them to $(x,x,0)$.
\end{itemize}

We shall next construct the main graph on the vertex set $X$. Recall that each directed labeled edge in the Cayley graph of $\group$ with respect to $S$, which corresponds to an element of $S$, is to be replaced by one of the continuum-many asymmetric graphs that we initially constructed. This supply of asymmetric graphs were parametrized by $\cantor$. Consequently, it suffices to parametrize $\group$ by $\cantor$. Since $(\group,\mathcal{B})$ is an uncountable standard Borel space, it follows from the Borel isomorphism theorem \cite[Theorem 15.6]{Kechris95} that there exists a Borel isomorphism $\Psi: \group \rightarrow \cantor$.

Before we proceed, we would like to take a moment to let the reader know in advance that we will later require $\Psi: \group \rightarrow \cantor$ to have other additional properties in the proof of Theorem \ref{maintheorem2}. Indeed, as we shall see later, the Borel complexity of our graph, i.e. where it resides in the Borel hierarchy of the Polish space $(X \times X,\tau \times \tau)$, is completely determined by the Borel complexity of inverse images of the clopen basis elements of $\cantor$ under $\Psi$.

\begin{figure}
  \includegraphics[width=\textwidth,height=\textheight,keepaspectratio]{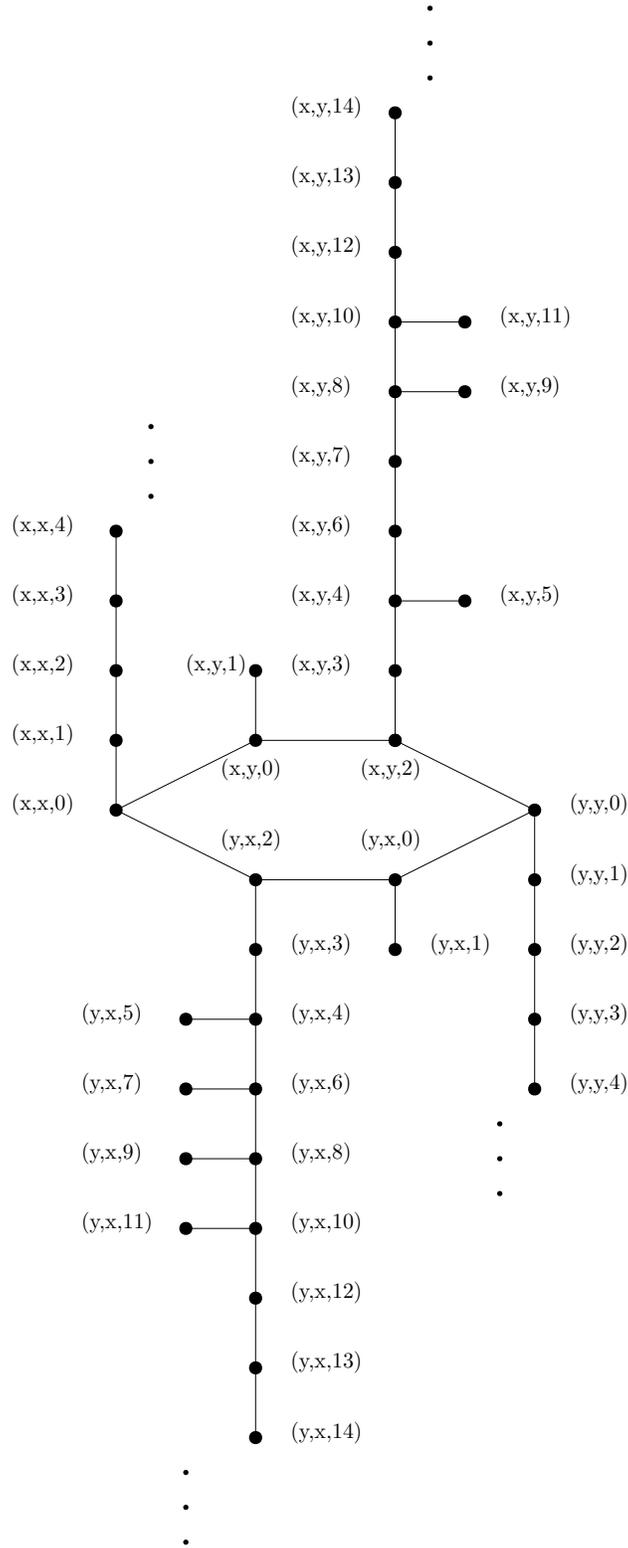}
  \caption{A representation of edges in $G^*$ for a pair of group elements $x$ and $y$}
  \label{figuremain}
\end{figure}

Consider the relation $G=G_{\text{irrelevant}}\cup  G_{\text{blockbase}}\cup  G_{\text{fork}}\cup G_{\text{nofork}}$ where\\
\begin{align*}
    G_{\text{blockbase}}=&\left\{\bigg( (x,x,0),(x,y,0)\bigg),\bigg((x,y,0),(x,y,1)\bigg),\bigg((x,y,0),(x,y,2)\bigg),\right.\\&\left.\ \ \bigg((x,y,2),(x,y,3)\bigg),\bigg((x,y,3),(x,y,4)\bigg),\bigg((x,y,2),(y,y,0)\bigg): x \neq y\in\group\right\}\\
    &\\
    G_{\text{forks}}=&\left\{\bigg((x,y,n),(x,y,n+1)\bigg),\bigg((x,y,n),(x,y,n+2)\bigg):\right.\\&\left.\ \ \ x \neq y\in\group,\ n \in 2\mathbb{N}_{\geq 2},\ \Psi\left(x^{-1}y\right)\left(\frac{n-4}{2}\right)=1\right\}\\
    &\\
    G_{\text{noforks}}=&\left\{\bigg((x,y,n),(x,y,n+1)\bigg),\bigg((x,y,n+1),(x,y,n+2)\bigg):\right.\\&\left.\ \ \ x \neq y\in \group,\ n \in 2\mathbb{N}_{\geq 2},\ \Psi\left(x^{-1}y\right)\left(\frac{n-4}{2}\right)=0\right\}\\
    &\\
    G_{\text{irrelevant}}=&\left\{\bigg((x,x,n),(x,x,n+1)\bigg):\ x\in \group,\ n\in\mathbb{N}\right\}\\
\end{align*}
An illustration of the edges in $G^*$ for a pair of group elements $x$ and $y$ is given in Figure \ref{figuremain} as an undirected graph, where we assume for illustrative purposes that $\Psi\left(x^{-1}y\right)=(1,0,1,1,0,\dots)$ and $\Psi\left(y^{-1}x\right)=(1,1,1,1,0,\dots)$.

It is a routine verification to check that $G$ is a Borel subset of $X \times X$. For those who are not well-versed in such matters, we will show here that $G_{\text{forks}}$ is indeed Borel as a guiding example. Let $n \in 2\mathbb{N}_{\geq 2}$. Then the set
\[O=\left\{\mathbf{a} \in \cantor: \mathbf{a}\left(\frac{n-4}{2}\right)=1\right\}\]
is a clopen subset of $\cantor=\{0,1\}^{\mathbb{N}}$. Consider the map from $f: \group \times \group \rightarrow \cantor$ given by $f(x,y)=\Psi\left(x^{-1}y\right)$. Since $(\group,\cdot,\mathcal{B})$ is a standard Borel group, $f$ is a Borel map and hence $B=f^{-1}(O)$ is a Borel subset of $\group \times \group$. Set $A=B \setminus \Delta_{\group}$. It follows that $A \times \{n\}$, $A \times \{n+1\}$ and $A \times \{n+2\}$ are Borel subsets of $X$ and hence, their pairwise cartesian products are Borel subsets of $X \times X$. But then
\[G_{\text{forks}}=\bigcup_{n \in 2\mathbb{N}_{\geq 2}} \big((A \times \{n\}) \times (A \times \{n+1\})\big) \cup \big((A \times \{n\}) \times (A \times \{n+2\})\big)\]
is a Borel subset of $X \times X$. That $G_{\text{noforks}}$, $G_{\text{blockbase}}$ and $G_{\text{irrelevant}}$ are Borel can be shown by similar arguments with appropriate modifications. Thus $\graph=(X,G^*)$ is a Borel graph.

\section{Proof of Theorem \ref{maintheorem}}

Let $(\group,\cdot,\mathcal{B})$ be a standard Borel group. Suppose for the moment that $\group$ is uncountable. Set $\graph=(X,G^*)$ to be the Borel graph constructed in Section 2 associated to $(\group,\cdot,\mathcal{B})$. We wish to show that $\group$ and $Aut_B(\graph)$ are isomorphic.  For each $g \in \group$, consider the map $\varphi_g: X \rightarrow X$ given by
\[\varphi_g(x,y,k)=(gx,gy,k)\]
for all $x,y \in \group$ and $k \in \mathbb{N}$. Clearly $\varphi_g$ is a bijective map. Observe that left-multiplying the first two components of each element of $G$ by $g$ leaves the sets $G_{\text{irrelevant}}$, $G_{\text{blockbase}}$, $G_{\text{fork}}$ and $G_{\text{nofork}}$ invariant. To see that $G_{\text{fork}}$ and $G_{\text{nofork}}$ are invariant under $\varphi_g$, observe that $x^{-1}y=(gx)^{-1}(gy)$. Thus $\varphi_g$ is an automorphism. Since the group multiplication is Borel, so is $\varphi_g$. It follows that $\varphi_g \in Aut_B(\graph)$. Define the map $\Phi: \group \rightarrow Aut_B(\graph)$ by $\Phi(g)=\varphi_{g}$ for all $g \in \group$. Then clearly $\Phi$ is injective and moreover, we have \[\Phi(gh)=\varphi_{gh}=\varphi_{g} \circ \varphi_{h}=\Phi(g) \circ \Phi(h)\]
Thus $\Phi$ is a group embedding. It remains to show that $\Phi$ is surjective.

Let $f \in Aut(\graph)$ be an arbitrary automorphism. Observe that the set of vertices which has one neighbor of infinite degree and another neighbor of degree one is precisely
\[\{(x,y,0): x \neq y \in \group\}\]
Therefore, being an automorphism, $f$ permutes this set. Let $x, y \in \group$ be distinct and set $(x',y',0)=f(x,y,0)$. Note that the only neighbors of $(x,y,0)$ and $(x',y',0)$
\begin{itemize}
    \item of degree one is $(x,y,1)$ and $(x',y',1)$ respectively,
    \item of degree three is $(x,y,2)$ and $(x',y',2)$ respectively,
    \item of uncountable degree is $(x,x,0)$ and $(x',x',0)$ respectively.
\end{itemize}
Thus $f(x,y,1)=(x',y',1)$, $f(x,y,2)=(x',y',2)$ and $f(x,x,0)=(x',x',0)$. By a similar argument, since we already obtained $f(x,y,2)=(x',y',2)$, we must also have that $f(y,y,0)=(y',y',0)$.

Recall that the graph $\graph_{\Psi(x^{-1}y)}=\left(\mathbb{N}_{\geq 2},R^*_{\Psi(x^{-1}y)}\right)$ constructed at the very beginning has no non-trivial automorphisms. Consequently, an inductive argument as was done in Section 2 shows that $f(x,y,n)=(x',y',n)$ for all $n \geq 2$.

This last conclusion immediately implies that $\Psi(x^{-1}y)=\Psi(x'^{-1}y')$. Since $\Psi$ is injective, we have $x^{-1}y=x'^{-1}y'$ and hence $x'x^{-1}=y'y^{-1}$.

Set $g=x'x^{-1} \in \group$. Then we have $gx=x'$ and $gy=y'$. Therefore
\[f(x,y,n)=(x',y',n)=\varphi_g(x,y,n)\]
for all $n \in \mathbb{N}$. Observe that if we used another $z \in \group$ instead of $x$ or $y$, we still would have found the same group element $g$ because in this case we would have $g=x'x^{-1}=y'y^{-1}=z'z^{-1}$. Therefore, we indeed have
\[f(x,y,n)=(x',y',n)=\varphi_g(x,y,n)\]
not only for the previously fixed $x,y$ but for all distinct $x,y \in \group$ and $n \in \mathbb{N}$. Hence $f$ agrees with $\varphi_g$ on $X \setminus (\Delta_{\group} \times \mathbb{N})$. It also follows from $f(x,x,0)=(gx,gx,0)$ via an inductive argument that $f(x,x,n)=(gx,gx,n)$ for all $n \in \mathbb{N}$. Thus $f$ is identically $\varphi_g$ on $X$. Hence $\Phi$ is an isomorphism and we indeed have $Aut(\graph)=Aut_B(\graph)$.

Finally, suppose that $\group$ is a countable standard Borel group. In this case, we choose $\Psi$ to be any (necessarily Borel) bijection from $\group$ to any (necessarily Borel) subset of $\cantor$ with cardinality $|\group|$ and implement the same construction. The exact same argument proving $Aut(\graph)=Aut_B(\graph)$ in the uncountable case still goes through in the countable case, with appropriate modifications in the extreme case $|\group|=1$.

\section{Proof of Theorem 2}

Let $(\group,\tau)$ be a Polish group. It is well-known \cite[Theorem 4.14]{Kechris95} that there exists a continuous injection $\gamma: \group \rightarrow \hilbert$. Consider the map $\xi: \hilbert \rightarrow 2^{\mathbb{N} \times \mathbb{N}}$ given by $\xi\left(\mathbf{x}\right)(i,j)=1$ if and only if the $i$-th digit of the binary expansion of $x_j$ is equal to 1, where the binary expansions of dyadic rationals are taken to end in infinitely many repating $1$'s. It is straightforward to check that $\xi$ is a $\mathbf{\Delta^0_2}$-map i.e. the inverse images of open sets are $\mathbf{\Delta^0_2}$. Recall that the Cantor pairing map from $\mathbb{N} \times \mathbb{N}$ to $\mathbb{N}$ given by $(m,n) \mapsto \frac{1}{2}(m+n)(m+n+1)+n$ is a bijection. So the map $\zeta: 2^{\mathbb{N} \times \mathbb{N}} \rightarrow \cantor$ given by
\[ \zeta(\mathbf{a})\left(\frac{1}{2}(m+n)(m+n+1)+n\right)=\mathbf{a}(m,n)\]
is a homeomorphism. Set $\widehat{\Psi}=\zeta \circ \xi \circ \gamma$.

We now carry out the same construction of $\graph=(X,G^{*})$ in Section 2 but we use the continuous injection $\widehat{\Psi}: \group \rightarrow \cantor$ instead of the Borel bijection $\Psi: \group \rightarrow \cantor$. Then the set $A$ in the construction is $\mathbf{\Delta^0_2}$. It follows that $G_{\text{forks}}$ and $G_{\text{noforks}}$ are $\mathbf{\Delta^0_2}$. It is also easily seen that $G_{\text{blockbase}}$ and $G_{\text{irrelevant}}$ are closed sets. Therefore $G^{*}$ is a $\mathbf{\Delta^0_2}$-subset of $X \times X$.

We next execute the proof of Theorem \ref{maintheorem} as it is. Observe that the automorphisms $\varphi_g: X \rightarrow X$ constructed in the proof are homeomorphisms. Moreover, $\widehat{\Psi}$ being injective suffices for the argument to go through. Thus we obtain that $\Aut(\graph)=\Aut_h(\graph)$ and that $\Phi: \group \rightarrow \Aut_h(\graph)$ is an isomorphism.

We shall next prove that $\Phi$ is indeed a homeomorphism whenever the group $\Aut_h(\graph) \subseteq \Homeo(X)$ is endowed with the subspace topology induced from the compact-open topology of $\Homeo(X)$. Let $\{O_{\alpha}\}_{\alpha \in I} \subseteq X$ be a basis for the topology of $X$. Recall that the collection
\[\{\{f \in C(X,X):\ f[K] \subseteq O_{\alpha}\}:\ K \subseteq X \text{ is compact}\}\]
is a subbase for the compact-open topology of $C(X,X)$. Let $U\subseteq \group$ be open. Then
\[\Phi(U)=\left\{\varphi_g\in \Aut_h(\graph):\ \varphi_g(1_\group,1_\group,1)\in U\times U\times \{1\} \right\}\]
Since the set $\{(1_\group,1_\group,1)\}$ is compact and $U\times U\times \{1\}$ is open in $X$, the set $\Phi(U)$ is open in the subspace topology of $\Aut_h(\graph)$. Hence $\Phi^{-1}$ is continuous.

Let $V_{K,O} \subseteq X$ be a subbasis element of the subspace topology of $\Aut_h(\graph)$ where $K \subseteq X$ is compact, $U_1 \times U_2 \times U_3=O \subseteq X$ is a basis element with $U_1,U_2 \subseteq \group$ and $U_3 \subseteq \mathbb{N}$ open; and
\[ V_{K,O}=\{\varphi \in \Aut_h(\graph): \varphi[K] \subseteq O\}\]
We wish to show that $\Phi^{-1}[V_{K,0}]=\{g \in \group: \varphi_g[K] \subseteq U_1 \times U_2 \times U_3\}$ is open. Observe that if $\pi_3[K] \nsubseteq U_3$, then $V_{K,O}=\emptyset$. So suppose that $\pi_{3}[K] \subseteq U_3$. Then we have
\[\Phi^{-1}[V_{K,0}]=\{g \in \group:  g \pi_1[K] \subseteq U_1\} \cap \{g \in \group:  g \pi_2[K] \subseteq U_2\}\]
We claim that both sets on the right hand side are open. To see this, let $g \in \group$ be such that $g \pi_i[K] \subseteq U_i$. For each $k \in \pi_i[K]$, since the multiplication on $\group$ is continuous and $gk \in U_i$, we can choose an open basis element $(g,k) \in V_k \times W_k$ of $\group \times \group$ such that $V_k \cdot W_k \subseteq U_i$. Since $\{W_k\}_{k \in K}$ is an open cover of the compact set $\pi_i[K]$, there exists a finite subcover $\{W_{k_i}\}_{i=1}^n$. Set $V=\bigcap_{i=1}^n V_{k_i}$. Then $g \in V$ and $V \cdot \pi_{i}[K] \subseteq U_i$. Thus $\{g \in \group:  g \pi_i[K] \subseteq U_i\}$ is open. Hence $\Phi$ is continuous and so, is a homeomorphism.

\section{Conclusion and further questions}

In this paper, we provided a complete generalization of Frucht's theorem to Borel measurable and topological settings. However, due to the natural limitations of our coding technique, in topological setting, we were not able to obtain minimal complexity in Theorem \ref{maintheorem2}. Therefore, we pose the following question.\\

\textbf{Question.} Is it true that for every Polish space $(\group,\tau)$ there exists a closed or open graph $\graph=(X,G)$ on a Polish space $(X,\widehat{\tau})$ such that $\group$ and $\Aut_h(\graph)$ are isomorphic?\\

We strongly suspect that the answer is affirmative. Such a result may be obtained via a construction similar to ours that uses a continuous injection $\Psi: \group \rightarrow [0,1]^{\mathbb{N}}$ which we know exists for arbitrary second-countable metrizable spaces $\group$. However, this would require one to construct continuum-many acyclic Borel graphs that code each element of the Hilbert cube $[0,1]^{\mathbb{N}}$ in such a way that each edge corresponds to an open or closed condition in $[0,1]^{\mathbb{N}}$. It is not clear to us how this can be done.

Observe that, due to the nature of the construction, the graphs that we obtained automatically ended up satisfying $\Aut(\graph)=\Aut_B(\graph)$. However, it is trivial to observe via counting arguments that it is possible to have Borel graphs $\graph$ such that $|\Aut_{B}(\graph)| \leq 2^{\aleph_0} < 2^{2^{\aleph_0}} \leq |\Aut(\graph)|$, for example, consider the complete graph $K_{\mathbb{R}}$. The next obvious question would be to ask whether it is possible to have $|\Aut_{B}(\graph)| \leq \aleph_0 < 2^{\aleph_0} \leq |\Aut(\graph)|$ for a Borel graph. Having corresponded with Andrew Marks, we learned that this question also has an affirmative answer. Here we briefly sketch his argument: Given a countable language $\mathcal{L}$, for any $\mathcal{L}$-structure $\mathcal{M}$ whose universe is a Polish space and whose functions and relations are Borel maps, one can construct a Borel graph $\graph_{\mathcal{M}}$ on a Polish space such that $\Aut_B(\mathcal{M}) \cong \Aut_B(\graph_{\mathcal{M}})$ and $\Aut(\mathcal{M}) \cong \Aut(\graph_{\mathcal{M}})$. This can be achieved by appropriately modifying the argument which shows that arbitrary structures may be interpreted as graphs, e.g. see \cite[Theorem 5.5.1]{Hodges93}. Consequently, it suffices to find $\mathcal{M}$ such that $|\Aut_{B}(\mathcal{M})| \leq \aleph_0 < 2^{\aleph_0} \leq |\Aut(\mathcal{M})|$.

An example of such a structure would be $(\mathbb{R},+,1)$. Since any Borel measurable group automorphism of the Polish group $(\mathbb{R},+)$ is automatically continuous \cite[Theorem 2.2]{Rosendal09} and any continuous automorphism of $(\mathbb{R},+)$ is precisely of the form $x \mapsto rx$, we have that $|\Aut_B(\mathbb{R},+,1)|=1$. On the other hand, since any permutation of a $\mathbb{Q}$-basis of $\mathbb{R}$ would induce a group automorphism and $\text{dim}_{\mathbb{Q}}(\mathbb{R})=2^{\aleph_0}$, we have $|\Aut(\mathbb{R},+,1)|=2^{2^{\aleph_0}}$.

Having seen that the Borel and full automorphism groups of a Borel graph can be separated in cardinality, the following question seems to be the next step in our initial investigation.

\textbf{Question.} Given two standard Borel groups $\mathcal{H} \leq \group$, does there necessarily exist a Borel graph $\graph$ such that $\Aut_B(\graph) \cong \mathcal{H}$ and $\Aut(\graph) \cong \group$, where the former isomorphism is the restriction of the latter?

\bibliography{references}{}

\begin{thebibliography}{KST99}

\bibitem[Bil22]{Bilge22}
Onur Bilge.
\newblock Some results on automorphism groups and dichotomy theorems in
  descriptive graph combinatorics.
\newblock Master's thesis, Middle East Technical University, 2022.
\newblock expected to be completed in July 2022.

\bibitem[dG59]{Groot59}
J.~de~Groot.
\newblock Groups represented by homeomorphism groups.
\newblock {\em Math. Ann.}, 138:80--102, 1959.

\bibitem[Fru39]{Frucht39}
R.~Frucht.
\newblock Herstellung von {G}raphen mit vorgegebener abstrakter {G}ruppe.
\newblock {\em Compositio Math.}, 6:239--250, 1939.

\bibitem[Hod93]{Hodges93}
Wilfrid Hodges.
\newblock {\em Model theory}, volume~42 of {\em Encyclopedia of Mathematics and
  its Applications}.
\newblock Cambridge University Press, Cambridge, 1993.

\bibitem[Kec95]{Kechris95}
Alexander~S. Kechris.
\newblock {\em Classical descriptive set theory}, volume 156 of {\em Graduate
  Texts in Mathematics}.
\newblock Springer-Verlag, New York, 1995.

\bibitem[KM20]{KechrisMarks20}
Alexander~S. Kechris and Andrew~S. Marks.
\newblock Descriptive graph combinatorics.
\newblock
  \url{http://www.math.caltech.edu/~kechris/papers/combinatorics20book.pdf},
  October 2020.

\bibitem[KST99]{KST99}
A.~S. Kechris, S.~Solecki, and S.~Todorcevic.
\newblock Borel chromatic numbers.
\newblock {\em Adv. Math.}, 141(1):1--44, 1999.

\bibitem[Ros09]{Rosendal09}
Christian Rosendal.
\newblock Automatic continuity of group homomorphisms.
\newblock {\em Bull. Symbolic Logic}, 15(2):184--214, 2009.

\bibitem[Sab60]{Sabidussi60}
Gert Sabidussi.
\newblock Graphs with given infinite group.
\newblock {\em Monatsh. Math.}, 64:64--67, 1960.

\end{thebibliography}
\bibliographystyle{alpha}

\end{document}